\documentclass[12pt,a4paper,leqno]{article}

\usepackage{theorem} \usepackage{graphics} \usepackage{amsfonts}

\theoremstyle{change} {\theorembodyfont{\itshape}
  \newtheorem{theorem}{Theorem.}[section]
  \newtheorem{proposition}[theorem]{Proposition.}
  \newtheorem{lemma}[theorem]{Lemma.}
   }

{\theorembodyfont{\rmfamily} \newtheorem{remark}[theorem]{Remark.}
  \newtheorem{example}[theorem]{Example.} }

\def\proof{\noindent{\it Proof.}\enspace} \def\endproof{ \quad
  $\kasten$}

\def\tq{\mathop{\lower 4pt \hbox{$ {\buildrel{/} \over
        {\scriptstyle{\rm tor}}} $}}\nolimits} \def\Xqp{X^{\rm qp}}
\def\Xtqp{X^{\rm tqp}}  \def\t#1{\widetilde{#1}}
\def\b#1{\overline{#1}} \def\CC{{\mathbb C}} \def\ZZ{{\mathbb Z}}
\def\RR{{\mathbb R}} \def\PP{{\mathbb P}} \def\mal{\mathbin{\! \cdot
    \!}}  
\def\conv{\mathop{\hbox{\rm conv}}} \def\cone{\mathop{\hbox{\rm
      cone}}} \def\id{\mathop{\rm id}\nolimits} \def\Hom{\mathop{\rm
    Hom}\nolimits}

\def\kasten{\mathord{\vbox{\hrule \hbox{\vrule \hskip5pt \vrule
        height5pt \vrule} \hrule}}}

\begin{document}

\thispagestyle{empty}

\begin{center}
  
  {\LARGE\bf Quasi-Projective Reduction}
  
  \medskip
  
  {\LARGE\bf of Toric Varieties}
  
  \bigskip
  
  Annette A'Campo--Neuen and J\"urgen Hausen
  
  \bigskip

\end{center}

\begin{abstract}
\noindent We define a quasi--projective reduction of a complex algebraic
  variety $X$ to be a regular map from $X$ to a quasi--projective
  variety that is universal with respect to regular maps from $X$ to
  quasi--projective varieties. A toric quasi--projective reduction is the
  analogous notion in the category of toric varieties. For a given
  toric variety $X$ we first construct a toric quasi--projective
  reduction. Then we show that $X$ has a quasi--projective reduction if
  and only if its toric quasi--projective reduction is surjective. We
  apply this result to characterize when the action of a subtorus on a
  quasi--projective toric variety admits a categorical quotient in the
  category of quasi--projective varieties.
\end{abstract}

\section*{Introduction}

Let $X$ be a complex algebraic variety. We call a regular map $p\colon X
\to \Xqp$ from $X$ to a quasi--projective variety $\Xqp$ a {\it
  quasi--projective reduction} of $X$ if every regular map $f \colon X \to
Z$ to a quasi--projective variety $Z$ factors uniquely through $p$,
i.e., there exists a unique regular map $\t{f} \colon \Xqp \to Z$ such $f =
\t{f} \circ p$.

As a first result of the present article we characterize when a given
toric variety $X$ has a quasi--projective reduction. In Section 1 we
construct a {\it toric quasi--projective reduction} of $X$, i.e., a
toric morphism $q\colon X \to \Xtqp$ to a quasi--projective toric variety
$\Xtqp$ such that every toric morphism from $X$ to a quasi--projective
toric variety factors uniquely through $q$. Then we prove (see Section
2):
 
\bigskip
 
\noindent{\bf Theorem 1.}\enspace {\it A toric variety $X$ has a
  quasi--projective reduction if and only if its toric quasi--projective
  reduction $q\colon X \to \Xtqp$ is surjective. If $q$ is surjective, then
  it is the quasi--projective reduction of $X$.}

\bigskip

The above theorem implies in particular that every complete toric
variety has a projective reduction. But as we show by an explicit
example (see \ref{nonsurj}), the quasi--projective reduction need not
exist in general. We apply Theorem 1 to obtain a complete answer to
the following problem, posed by A. Bia\l ynicki-Birula:

Let $X$ be a quasi--projective toric variety with acting torus $T$ and
let $H \subset T$ be a subtorus. When does the action of $H$ admit a
quotient in the category of quasi--projective varieties, i.e., an
$H$-invariant regular map $s\colon X \to Y$ to a quasi--projective variety
$Y$ such that every $H$-invariant regular map from $X$ to a
quasi--projective variety factors uniquely through $s$?

In order to state our answer, let $s_{1} \colon X \to X \tq H$ denote the
toric quotient (see \cite{acha}). Recall that $s_{1}$ is universal
with respect to $H$-invariant toric morphisms. Moreover, let $q\colon X \tq
H \to Y$ be the toric quasi--projective reduction. Then our result is
the following (for the proof see Section 2):
 
\goodbreak\bigskip
 
\noindent{\bf Theorem 2.}\enspace {\it The action of $H$ on $X$ admits a
  quotient in the category of quasi--projective varieties if and only
  if $s := q \circ s_{1}$ is surjective. If $s$ is surjective, then it
  is the quotient for the action of $H$ on $X$.}

\bigskip

Examples of quasi--projective toric varieties with a subtorus action
admitting a quotient in the category of quasi--projective varieties are
obtained from Mumford's Geometric Invariant Theory. In \ref{noquot}
and \ref{subtil} we discuss examples of subtorus actions that have no
such quotient.

\section*{Notation}

A {\it toric variety\/} is a normal algebraic variety $X$ endowed with
an effective regular action of an algebraic torus $T$ that has an open
orbit. We refer to $T$ as the {\it acting torus\/} of $X$. For every
toric variety $X$ we fix a base point $x_{0}$ in the open orbit.

A regular map $f \colon X \to X'$ of toric varieties with base points
$x_{0}$ and $x_{0}'$ respectively is called a {\it toric morphism\/}
if $f(x_{0}) = f(x'_{0})$ and there is a homomorphism $\varphi \colon
T \to T'$ of the acting tori such that $f(t \mal x) = \varphi(t) \mal
f(x)$ holds for every $(t,x) \in T \times X$.

The basic construction in the theory of toric varieties is to
associate to a given fan $\Delta$ in an $n$-dimensional lattice an
$n$-dimensional toric variety $X_{\Delta}$. The assignment $\Delta
\mapsto X_{\Delta}$ is in fact an equivalence of categories (see e.g.
\cite{Fu} or \cite{Od}). For our construction of toric quasiprojective
reductions we need the following generalization of the notion of a
fan:

Let $N$ denote a $n$-dimensional lattice and set $N_{\RR} := \RR
\otimes_{\ZZ} N$. A {\it quasi--fan} in $N$ is a finite set $\Delta$ of
rational convex polyhedral cones in $N_{\RR}$ such that for each
$\sigma \in \Delta$ also every face of $\sigma$ is an element of
$\Delta$ and any two cones of $\Delta$ intersect in a common face. So
a quasi--fan is a fan if all its cones are strictly convex. For a
quasi--fan $\Delta$, we denote by $\Delta^{\max}$ the set of its
maximal cones and by $\vert \Delta \vert := \bigcup_{\sigma \in
  \Delta} \sigma$ its support.

For a homomorphism $F \colon N \to N'$ of lattices, let $F_{\RR}$
denote the associated homomorphism of real vector spaces. A {\it map
  of quasi--fans} $\Delta$ in $N$ and $\Delta'$ in $N'$ is a lattice
homomorphism $F \colon N \to N'$ such that for every $\sigma \in
\Delta$ there is a $\sigma' \in \Delta'$ with $F_{\RR}(\sigma) \subset
\sigma'$. As mentioned above, every map $F \colon \Delta \to \Delta'$
of fans gives rise to a toric morphism $f \colon X_{\Delta} \to
X_{\Delta'}$.

Every quasi--fan $\Delta$ in $N$ defines in a canonical manner a fan:
Let $V$ denote the intersection of all cones of $\Delta$. Then $V$ is
a linear subspace of $N_{\RR}$. Set $L := V \cap N$ and let $Q\colon N
\to \t{N} := N/L$ denote the projection. Then the cones
$Q_{\RR}(\sigma)$, $\sigma \in \Delta^{\max}$, are the maximal cones
of a fan $\t{\Delta}$ in $\t{N}$. We call $\t{\Delta}$ the {\it
  quotient fan} of $\Delta$.  By construction, $Q$ is a map of the
quasi--fans $\Delta$ and $\t{\Delta}$.

\section{Construction of the Toric Quasi-Projective \\ Reduction}

The construction of the toric quasi--projective reduction is done in
the category of fans. Toric morphisms from a complete toric variety
$X_{\Delta}$ to projective spaces are related to concave support
functions of the fan $\Delta$. Since we also want to consider
non-complete fans it is more natural to work with the following notion
instead of support functions:

\goodbreak

Let $\Delta$ be a quasi--fan in a lattice $N$. A finite family
${\mathfrak U} := (u_{i})_{i \in I}$ of linear forms $u_{i} \in M :=
\Hom(N,\ZZ)$ is called {\it $\Delta$-concave}, if it satisfies the
following condition: for every $\sigma \in \Delta^{\max}$ there is an
index $i(\sigma)$ such that
$$u_{i(\sigma)}\lower 4pt \hbox{$\vert \sigma$} \; \le \; u_{i}\lower
4pt \hbox{$\vert \sigma$} \quad \mbox{for all } i \in I.$$
Note that
for two given $\Delta$-concave families ${\mathfrak U} := (u_{i})_{i
  \in I}$ and ${\mathfrak U}' := (u'_{j})_{j \in J}$ of linear forms
the {\it sum family}
$$
{\mathfrak U} + {\mathfrak U}' := (u_{i} + u'_{j})_{(i,j) \in I
  \times J}$$
is again a $\Delta$-concave family. For a
$\Delta$-concave family ${\mathfrak U}$ let $P_{{\mathfrak U}}$ denote
the convex hull of ${\mathfrak U}$.  Then $P_{{\mathfrak U}}$ is a
lattice polytope in $M_{\RR}$. Let $\Sigma_{\mathfrak U}$ denote the
normal quasi--fan of $P_{\mathfrak U}$ in $N$. Recall that the faces
$P'$ of $P$ correspond order-reversingly to the cones of
$\Sigma_{\mathfrak U}$ by
$$
P' \mapsto \tau_{P'} := \{v \in N_{\RR}; \; p'(v) \le p(v) \hbox{
  for all } (p',p) \in P' \times P\}.$$
If $u_{1}, \ldots, u_{r}$
denote the vertices of $P_{\mathfrak U}$, then the family $(u_{i})_{i
  = 1, \ldots, r}$ is {\it strictly $\Sigma_{\mathfrak U}$-concave},
i.e., on every relative interior $\tau_{\{u_{i}\}}^{\circ}$ the linear
form $u_{i}$ is strictly smaller than the forms $u_{j}$ with $j \ne
i$.

Now assume that $\Delta$ is a fan. Call a subset $R$ of the set
$\Delta^{(1)}$ of extremal rays of $\Delta$ {\it indecomposable}, if
for every $\Delta$-concave family ${\mathfrak U}$ the set $R$ is
contained in some maximal cone of $\Sigma_{\mathfrak U}$. Let $R_{1},
\ldots, R_{k}$ be the maximal indecomposable subsets of
$\Delta^{(1)}$.

\begin{lemma}\label{generic}
  There exists a $\Delta$-concave family ${\mathfrak U}$ such that
  every $R_{i}$ is the intersection of $\Delta^{(1)}$ with some
  maximal cone of $\Sigma_{\mathfrak U}$.
\end{lemma}

\proof For every decomposable subset $S$ of $\Delta^{(1)}$ choose a
$\Delta$-concave family ${\mathfrak U}_{S}$ such that $S$ is not
contained in any maximal cone of $\Sigma_{{\mathfrak U}_{S}}$. Let
${\mathfrak U}$ be the sum of these families ${\mathfrak U}_{S}$. Then
$P_{\mathfrak U}$ is the Minkowski-Sum of the $P_{{\mathfrak U}_{S}}$.
Consequently $\Sigma_{\mathfrak U}$ is the common refinement of the
$\Sigma_{{\mathfrak U}_{S}}$. This readily yields the claim. \endproof

\bigskip

A $\Delta$-concave family ${\mathfrak U}$ with the property of Lemma
\ref{generic} will be called {\it generic}. As a consequence of the
above lemma we obtain the following statement for the sets
$$
\varrho_{i} := \conv\biggl(\bigcup_{\varrho \in R_{i}} \varrho
\biggr). $$

\begin{remark}\label{Sigma}
  The $\varrho_{i}$ are the maximal cones of a quasi--fan $\Sigma$ in
  $N$. The lattice homomorphism $\id_{N}$ is a map of the quasi--fans
  $\Delta$ and $\Sigma$. Moreover, if ${\mathfrak U}$ is a generic
  $\Delta$-concave family, then $\id_{N}$ is an affine map of the
  quasi--fans $\Sigma$ and $\Sigma_{\mathfrak U}$.  \endproof
\end{remark} 

Here a map $F$ of quasi--fans $\Delta$ in $N$ and $\Delta'$ in $N'$ is
called {\it affine} if for every maximal cone $\sigma'$ of $\Delta'$
the set $F_{\RR}^{-1}(\sigma') \cap \vert \Delta \vert$ is a (maximal)
cone of $\Delta$. Note that a map of fans is affine if and only if the
associated toric morphism is affine.

Now we construct the quasi--projective toric reduction of a toric
variety $X_{\Delta}$ defined by the fan $\Delta$. Let $V$ denote the
minimal cone of the quasi--fan $\Sigma$ determined by $\Delta$ as in
Remark \ref{Sigma}. Set $L := N \cap V$, let $Q \colon N \to \t{N} :=
N/L$ be the projection and denote by $\t{\Delta}$ the quotient-fan of
$\Sigma$.

\begin{proposition}
  The toric morphism $q \colon X_{\Delta} \to X_{\t{\Delta}}$
  associated to $Q$ is the toric quasi--projective reduction of
  $X_{\Delta}$.
\end{proposition}

\proof First we show that $X_{\t{\Delta}}$ is in fact
quasi--projective.  Choose a generic $\Delta$-concave family
${\mathfrak U} = (u_{\sigma})_{\sigma \in \Delta^{\max}}$. Let $V_1$
denote the minimal cone of the quasi--fan $\Sigma_{\mathfrak U}$. Set
$L_1 := N \cap V_1$, let $P \colon N \to \b{N} := N/L_{1}$ be the
projection and denote the quotient-fan of $\Sigma_{\mathfrak U}$ by
$\b{\Delta}$.

Since ${\mathfrak U}$ induces a strictly $\b{\Delta}$-concave family,
the associated toric variety $X_{\b{\Delta}}$ is projective.  The
minimal cone $V$ of $\Sigma$ is contained in $V_1$, so we obtain a
lattice homomorphism $G \colon \t{N} \to \b{N}$ with $G \circ Q = P$.
By construction, $G$ is an affine map of the fans $\t{\Delta}$ and
$\b{\Delta}$. So the associated toric morphism $g \colon
X_{\t{\Delta}} \to X_{\b{\Delta}}$ is affine. Since $X_{\b{\Delta}}$
is projective we can use \cite{EGA}, Chap. II, Th. 4.5.2, to conclude
that $X_{\t{\Delta}}$ is quasi--projective.

Now we verify the universal property of $q$. Let $f \colon X_{\Delta}
\to X'$ be a toric morphism to a quasi--projective toric variety $X'$.
We may assume that $f$ arises from a map $F \colon N \to N'$ of fans
$\Delta$ and $\Delta'$. Choose a polytopal completion $\Delta''$ of
$\Delta'$.  By suitable stellar subdivisions (see \cite{Ew}, p. 72) we
achieve that every maximal cone of $\Delta''$ contains at most one
maximal cone of $\Delta'$.

Let $(u_{\sigma''})_{\sigma'' \in {\Delta''}^{\max}}$ be a strictly
$\Delta''$-concave family. Then the linear forms $u_{\sigma''} \circ
F$ form a $\Delta$-concave family. Let $\sigma \in \Sigma^{\max}$.  By
construction, $\sigma$ is mapped by $F_{\RR}$ into some cone of
$\Delta''$. Moreover, $\sigma$ is the convex hull of certain extremal
rays of $\Delta$, so $F_{\RR}(\sigma)$ is in fact contained in a
maximal cone of $\Delta'$. Hence $F$ is a map of the quasi--fans
$\Sigma$ and $\Delta'$. In particular we have $F(L) = 0$. Thus there
is a map $\t{F} \colon \t{N} \to N'$ of the fans $\t{\Delta}$ and
$\Delta'$ with $F = \t{F} \circ Q$. The associated toric morphism
$\t{f} \colon X_{\t{\Delta}} \to X'$ yields the desired factorization
of $f$. \endproof

\section{Proof of the Theorems}

Let $X$ be a toric variety with acting torus $T$ and assume that $H
\subset T$ is an algebraic subgroup. Let $Z$ be an arbitrary
quasi--projective variety. We need the following decomposition result
for regular maps:

\begin{proposition}\label{factmor}
  Let $f \colon X \to Z$ be an $H$-invariant regular map. Then there
  exist a locally closed subvariety $W$ of some $\PP_{r}$, an
  $H$-invariant toric morphism $g\colon X \to \PP_{r}$ with $g(X)
  \subset W$ and a regular map $h\colon W \to Z$ such that $f = h
  \circ g$.
\end{proposition}

\proof In a first step we consider the special case that $Z = \PP_{m}$
and $X$ is an open toric subvariety of some $\CC^{n}$. Then there are
polynomials $f_0, \ldots, f_m\in\CC[z_1, \ldots, z_n]$ having no
common zero in $X$ such that $f(z) = [f_0(z), \ldots, f_m(z)]$ holds
for every $z \in X$. Clearly we may assume that the $f_{i}$ have no
non-trivial common divisor.

Since $f$ is $H$-invariant, every $f_i/f_j$ is an $H$-invariant
rational function on $\CC^{n}$. Thus, using $1 \in \gcd(f_0(z),
\ldots, f_m(z))$ we can conclude that there is a character $\chi
\colon H \to \CC^*$ satisfying $f_i(h \mal x) = \chi(h) f_i(x)$ for
every $i$ and every $(h,x)$ in $H \times X$.

Now every $f_i$ is a sum of monomials $q_{i1}, \ldots, q_{ir_i}$. Note
that also each of the monomials $q_{ij}$ is homogeneous with respect
to the character $\chi$. Moreover, since the $f_i$ have no common zero
in $X$, neither have the $q_{ij}$. Set $r := \sum r_i$ and define a
toric morphism
$$
g \colon X \to \PP_{r}, \quad x \mapsto [q_{01}(x), \ldots,
q_{0r_0}(x), \ldots, q_{m1}(x), \ldots, q_{mr_m}(x)]. $$
Then $g$ is $H$-invariant. In order to define an open subset $W$ of
$\PP_{r}$ and a regular map $h \colon W \to \PP_{m}$ with the desired
properties, consider the linear forms
$$
L_{i} \colon [z_{01}, \ldots, z_{0r_{0}}, \ldots, z_{m1}, \ldots,
z_{mr_{m}}] \mapsto z_{i1} + \ldots + z_{ir_{i}}$$
on $\PP_{r}$. Set
$W := \PP_{r} \setminus V(\PP_{r}; L_{1}, \ldots, L_{m})$ and
$$h \colon W \to \PP_{m}, \quad [z] \mapsto [L_{0}(z), \ldots,
L_{m}(z)].$$
Since the $f_{i}$ have no common zero in $X$ we obtain
$g(X) \subset W$. Moreover, by construction we have $f = h \circ g$.
So the assertion is proved for the case that $Z = \PP_{m}$ and $X$ is
an open toric subvariety of some $\CC^{n}$.

In a second step assume that $Z$ is arbitrary but $X$ again is an open
toric subvariety of some $\CC^{n}$. Choose a locally closed embedding
$\imath \colon Z \to \PP_{m}$. By Step one we obtain a decomposition
of $f' := \imath \circ f$ as $f' = h' \circ g$, where $g \colon X \to
\PP_{r}$ is an $H$-invariant toric morphism such that $g(X)$ is
contained in an open subset $W'$ of $\PP_{r}$ and $h' \colon W' \to
\PP_{m}$ is regular.

Then $W := {h'}^{-1}(\imath(Z))$ is a locally closed subvariety of
$W'$. Moreover we have $g(X) \subset W$ and there is a unique regular
map $h \colon W \to Z$ with $h' = \imath \circ h$. It follows that $f
= h \circ g$ is the desired decomposition.

Finally, let also $X$ be arbitrary. As described in \cite{cox}, there
is an open toric subvariety $U$ of some $\CC^n$ and a surjective toric
morphism $p \colon U \to X$ such that $p$ is the good quotient of $U$
by some algebraic subgroup $H_{0}$ of $(\CC^*)^n$.  Consider $f' := f
\circ p$. Then $f'$ is invariant by the action of $H' := \pi^{-1}(H)$,
where $\pi$ denotes the homomorphism of the acting tori associated to
$p$.

By the first two steps we can decompose $f'$ as $f' = h \circ g'$ with
an $H'$-invariant toric morphism $g' \colon U \to \PP_{r}$ and a
regular map $h \colon W \to Z$, where $W \subset \PP_{r}$ is locally
closed with $g'(U) \subset W$. Since $g'$ is $H'$-invariant, it is
also invariant by the action of $H_{0}$. Thus there is a unique toric
morphism $g\colon X \to \PP_{r}$ such that $g' = g \circ p$. Since $p$
is surjective, $g$ is $H$-invariant and we have $f = h \circ g$ which
is the desired decomposition of $f$. \endproof

\bigskip

\noindent{\it Proof of Theorem 1.}\enspace Let us first assume that
the toric quasi--projective reduction $q\colon X \to \Xtqp$ is
surjective. Let $f\colon X\to Z$ be a regular map to a
quasi--projective variety $Z$. We have to show that $f$ factors
uniquely through $q$.

By Proposition \ref{factmor} there is a toric morphism $g\colon X \to
X'$ to a projective toric variety $X'$, and a rational map $h\colon X'
\to Z$ which is regular on $g(X)$ such that $f = h \circ g$. Now there
is a toric morphism $\t{g} \colon \Xtqp \to X'$ such that $g = \t{g}
\circ q$.  Since $q$ was assumed to be surjective, we have
$\t{g}(\Xtqp) = g(X)$, and hence $f$ factors through $q$.

Now suppose that $p\colon X \to \Xqp$ is a quasi--projective reduction.
Then clearly $p$ is surjective and $\Xqp$ is normal. Moreover, there
is an induced action of the torus $T$ on $\Xqp$ making $p$
equivariant. We claim that this action is regular:

Accoding to Proposition \ref{factmor} choose a toric morphism $g\colon
X \to X'$ to a projective toric variety $X'$, and a rational map
$h\colon X' \to \Xqp$ such that $g(X)$ is contained in the domain $W'$
of definition of $h$ and $p=h\circ g$. By the universal property of
the toric quasi--projective reduction $q \colon X \to \Xtqp$ there is a
toric morphism $\t{g} \colon \Xtqp \to X'$ such that $g = \t{g} \circ
q$.

Moreover, by the universal property of $p$, there is a regular map
$\alpha \colon \Xqp \to \Xtqp$ such that $q = \alpha \circ p$. Note
that $\alpha(\Xqp) \subset q(X)$ and $\t{g}(q(X)) \subset W'$. So,
using surjectivity of $p$ and equivariance of $p$ and $q$, we obtain
for a given pair $(t,y) \in T \times \Xqp$ the equality
$$t \mal y = h(\t{g}(t \mal \alpha(y))).$$
This implies regularity of
the induced $T$-action on $\Xqp$. It follows that $\Xqp$ is in fact a
toric variety and $p$ is a toric morphism. Thus we obtain a toric
morphism $\beta \colon \Xtqp \to \Xqp$ with $p = \beta \circ q$. By
uniqueness of the factorizations we obtain that $\alpha$ and $\beta$
are inverse to each other, i.e., $q$ is also a quasi--projective
reduction of $X$.  \endproof

\bigskip

\noindent{\it Proof of Theorem 2.}\enspace Suppose first that $s
\colon X \to Y$ is a quotient for the action of $H$ on $X$ in the
category of quasi--projective varieties. As above we see that $Y$ is a
toric variety and $s$ is a surjective toric morphism. The universal
property of the toric quotient $s_{1} \colon X \to X \tq H$ yields a
toric morphism $q \colon X \tq H \to Y$ such that $s = q \circ
s_{1}$. Clearly $q$ satisfies the universal property of the toric
quasi--projective reduction of $X \tq H$.

Now let $s_{1} \colon X \to X \tq H$ denote the toric quotient, $q
\colon X \tq H \to Y$ the toric quasi--projective reduction and assume
that $q \circ s_{1}$ is surjective. Let $f \colon X \to Z$ be an
$H$-invariant regular map to a quasi--projective variety. Choose a
decomposition $f = h \circ g$ as in Proposition \ref{factmor}. Then,
by the universal properties of toric quotient and quasi--projective
reduction there is a regular map $\t{g}$ with $g = \t{g} \circ q \circ
s_{1}$. Since $q \circ s_{1}$ is surjective, $h$ is defined on
$\t{g}(Y)$. Thus $f = (h \circ \t{g}) \circ (q \circ s_{1})$ is the
desired factorization of $f$.  \endproof

\bigskip

The above proofs yield in fact the following generalization of
Theorems 1 and 2: Let $X$ be any toric variety and let $H$ be an
algebraic subgroup of the acting torus $T$ of $X$. Call an
$H$-invariant regular map $p \colon X \to X^{\rm qp}_H$ to a
quasi--projective variety $X^{\rm qp}_H$ an {\it $H$-invariant
  quasi--projective reduction} if it is universal with respect to
$H$-invariant regular maps from $X$ to quasi--projective varieties.

Now write $H = \Gamma H^0 $ with a finite subgroup $\Gamma$ and a
subtorus $H^0$ of $T$. Let $g \colon X \to X'$ denote the geometric
quotient for the action of $\Gamma$ on $X'$. Then $g$ is a toric
morphism. Hence there is an induced action of $H^0$ on $X'$. Let $s_1
\colon X' \to X' \tq H^0$ be the toric quotient for this action and
let $q \colon X' \tq H^0 \to Y$ be the quasi--projective toric
reduction. Then we obtain:

\bigskip

\noindent{\bf Theorem 3.}\enspace{\it $X$ has an $H$-invariant
  quasi--projective  reduction if and only $q \circ s_1$ is
  surjective. If so, then $q \circ s_1 \circ g$ is the $H$-invariant
  quasi--projective reduction of $X$. \endproof}

\section{Examples}

We first give an example of a $3$-dimensional toric variety
$X_{\Delta}$ that admits no quasi--projective reduction.  This variety
is an open toric subvariety of the minimal example for a smooth
complete but non-projective toric variety presented in \cite{Od},
Section 2.3.

\begin{example}\label{nonsurj}
  Let $e_{1}$, $e_{2}$ and $e_{3}$ denote the canonical basis vectors
  of the lattice $\ZZ^{3}$. Consider the vectors
  $$\begin{array}{ll}
    v_{1} := -e_{1}, & \qquad v_{1}' := e_{2} + e_{3}, \\
    v_{2} := -e_{2}, & \qquad v_{2}' := e_{1} + e_{3}, \\
    v_{3} := -e_{3}, & \qquad v_{3}' := e_{1} + e_{2}. \\
\end{array}$$
Let $\Delta$ be the fan in $\ZZ^3$ with the maximal cones
$$
\tau_1:=\cone(v_1,v_{3}'), \quad \tau_2:=\cone(v_2,v_{1}')\quad
\hbox{and} \quad \tau_3:=\cone(v_3,v_{2}').$$
\begin{center}
  \begin{picture}(0,0)%
\includegraphics{delta.pstex}%
\end{picture}%
\setlength{\unitlength}{0.00035000in}%
\begingroup\makeatletter\ifx\SetFigFont\undefined
\def\x#1#2#3#4#5#6#7\relax{\def\x{#1#2#3#4#5#6}}%
\expandafter\x\fmtname xxxxxx\relax \def\y{splain}%
\ifx\x\y   
\gdef\SetFigFont#1#2#3{%
  \ifnum #1<17\tiny\else \ifnum #1<20\small\else
  \ifnum #1<24\normalsize\else \ifnum #1<29\large\else
  \ifnum #1<34\Large\else \ifnum #1<41\LARGE\else
     \huge\fi\fi\fi\fi\fi\fi
  \csname #3\endcsname}%
\else
\gdef\SetFigFont#1#2#3{\begingroup
  \count@#1\relax \ifnum 25<\count@\count@25\fi
  \def\x{\endgroup\@setsize\SetFigFont{#2pt}}%
  \expandafter\x
    \csname \romannumeral\the\count@ pt\expandafter\endcsname
    \csname @\romannumeral\the\count@ pt\endcsname
  \csname #3\endcsname}%
\fi
\fi\endgroup
\begin{picture}(4500,4050)(4141,-3796)
\put(5131,-3166){\makebox(0,0)[lb]{\smash{\SetFigFont{8}{9.6}{rm}$v_1'$}}}
\put(4141,-3751){\makebox(0,0)[lb]{\smash{\SetFigFont{8}{9.6}{rm}$v_1$}}}
\put(6346, 74){\makebox(0,0)[lb]{\smash{\SetFigFont{8}{9.6}{rm}$v_3$}}}
\put(6346,-1141){\makebox(0,0)[lb]{\smash{\SetFigFont{8}{9.6}{rm}$v_3'$}}}
\put(7696,-3166){\makebox(0,0)[lb]{\smash{\SetFigFont{8}{9.6}{rm}$v_2'$}}}
\put(8641,-3751){\makebox(0,0)[lb]{\smash{\SetFigFont{8}{9.6}{rm}$v_2$}}}
\end{picture}

\end{center}
We claim that the toric quasi--projective reduction $q$ of $X_{\Delta}$
is the toric morphism associated to $\id_{N}$ interpreted as a map
from $\Delta$ to the fan $\t{\Delta}$ having as its maximal cones
$$
\sigma_1:=\cone(v_1,v_3,v_{1}',v_{3}'), \quad
\sigma_2:=\cone(v_1,v_2,v_{1}',v_{2}')\quad\hbox{and}\quad
\sigma_3:=\cone(v_2,v_3,v_{2}',v_{3}').
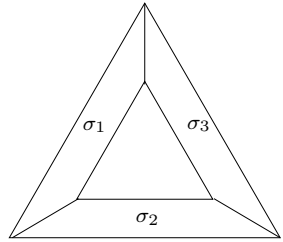
$$
\begin{center}
  \begin{picture}(0,0)%
\includegraphics{tdelta.pstex}%
\end{picture}%
\setlength{\unitlength}{0.00035000in}%
\begingroup\makeatletter\ifx\SetFigFont\undefined
\def\x#1#2#3#4#5#6#7\relax{\def\x{#1#2#3#4#5#6}}%
\expandafter\x\fmtname xxxxxx\relax \def\y{splain}%
\ifx\x\y   
\gdef\SetFigFont#1#2#3{%
  \ifnum #1<17\tiny\else \ifnum #1<20\small\else
  \ifnum #1<24\normalsize\else \ifnum #1<29\large\else
  \ifnum #1<34\Large\else \ifnum #1<41\LARGE\else
     \huge\fi\fi\fi\fi\fi\fi
  \csname #3\endcsname}%
\else
\gdef\SetFigFont#1#2#3{\begingroup
  \count@#1\relax \ifnum 25<\count@\count@25\fi
  \def\x{\endgroup\@setsize\SetFigFont{#2pt}}%
  \expandafter\x
    \csname \romannumeral\the\count@ pt\expandafter\endcsname
    \csname @\romannumeral\the\count@ pt\endcsname
  \csname #3\endcsname}%
\fi
\fi\endgroup
\begin{picture}(4094,3554)(4479,-3683)
\put(5626,-2041){\makebox(0,0)[lb]{\smash{\SetFigFont{8}{9.6}{rm}$\sigma_1$}}}
\put(7186,-2026){\makebox(0,0)[lb]{\smash{\SetFigFont{8}{9.6}{rm}$\sigma_3$}}}
\put(6421,-3421){\makebox(0,0)[lb]{\smash{\SetFigFont{8}{9.6}{rm}$\sigma_2$}}}
\end{picture}

\end{center}
Note that $q$ is not surjective. In order to prove that $q$ is the
toric quasi--projective reduction of $X_{\Delta}$, we have to show that
every $\Delta$-concave family $(u_{i})_{i=1,2,3}$ can be extended to a
$\t{\Delta}$-concave family. Note that $v_{1} + v_{3}'$ equals $v_{3}
+ v_{1}'$ and hence we have
$$\begin{array}{ccccc} u_{1}(v_{1}) + u_{1}(v_{3}') & = & u_{1}(v_{3})
  + u_{1}(v_{1}') & \ge &
  u_{3}(v_{3}) + u_{2}(v_{1}'). \\
\end{array}$$
Similarly we obtain
$$\begin{array}{ccccc} u_{2}(v_{2}) + u_{2}(v_{1}') & = & u_{2}(v_{1})
  + u_{2}(v_{2}') & \ge &
  u_{1}(v_{1}) + u_{3}(v_{2}'), \\
  u_{3}(v_{3}) + u_{3}(v_{2}') & = & u_{3}(v_{3}') + u_{3}(v_{2}) &
  \ge &
  u_{1}(v_{3}') + u_{2}(v_{2}). \\
\end{array}$$
Summing over these three inequalities, we arrive at an identity, and
therefore the inequalities are in fact equalities. This implies
$$\begin{array}{ll}
  u_{1}(v_{1}) = u_{2}(v_{1}), & \qquad u_{1}(v_{1}') = u_{2}(v_{1}'), \\
  u_{1}(v_{3}) = u_{3}(v_{3}), & \qquad u_{1}(v_{3}') = u_{3}(v_{3}'), \\
  u_{2}(v_{2}) = u_{3}(v_{2}), & \qquad u_{2}(v_{2}') = u_{3}(v_{2}').
  \quad \diamondsuit \\ \end{array}$$
\end{example}

In the above example the quasi--projective toric reduction has a
trivial kernel, and the variety $X_{\t{\Delta}}$ has the same
dimension as $X_{\Delta}$. For the complete case we have more
generally:

\begin{remark}
  Let $\Delta$ be a complete fan in a lattice $N$. Then $\dim
  X_{\Delta} = \dim X_{\Delta}^{\rm qp}$ holds if and only if $\Delta$
  can be defined via a subdivision of a lattice polytope in
  $N_{\RR}$.\endproof
\end{remark}

The next example is taken from the book of Fulton. It shows that in
general a complete toric variety is very far from its projective
reduction.

\begin{example} 
  Consider the complete fan in $\ZZ^{3}$ obtained by taking the cones
  over the faces of the standard cube with vertices $(\pm 1,\pm 1,\pm
  1)$. Deform this fan into a new complete fan $\Delta$ in $\ZZ^{3}$
  by moving the vertex $(1,1,1)$ to $(1,2,3)$.  The only support
  functions of $\Delta$ are the linear functions in $M$ (see
  \cite{Fu}, p. 26). So $X_{\Delta}^{\rm qp}$ is a point.  \quad
  $\diamondsuit$
\end{example}

Now we turn to quotients of a quasi--projective toric variety $X$ with
acting torus $T$ by subtori $H \subset T$. Examples of such quotients
are obtained by Mumford's Geometric Invariant Theory:

For the sake of simplicity assume $X=\PP_{n}$. Then the choice of a
lifting of the $T$-action to $\CC^{n+1}$ yields a notion of
$H$-semistability. The set $X^{\rm ss} \subset X$ of $H$-semistable
points is $T$-invariant and there is a quotient $X^{\rm ss} \to Y$ in
the category of quasi--projective varieties for the action of $H$ on
$X^{\rm ss}$ (see \cite{Mu}, also \cite{KaStZe} and \cite{BBSw}).

\begin{example}\label{noquot}
  Let $\Delta$ be the fan in $\RR^{4}$ that has $\sigma_{1} :=
  \cone(e_{1},e_{2})$ and $\sigma_{2} := \cone(e_{3},e_{4})$ as
  maximal cones. Then $X_{\Delta}$ is an open toric subvariety of
  $\CC^{4}$ with acting torus $T = {\CC^{*}}^{4}$. Define a projection
  $S_{1} \colon \ZZ^{4} \to \ZZ^{3}$ by setting
  $$S_{1}(e_{1}) := e_{1}, \quad S_{1}(e_{2}) := e_{2}, \quad
  S_{1}(e_{3}) := e_{3}, \quad S_{1}(e_{4}) := e_{1} + e_{2}.$$
  Then
  $S_{1}(e_{1}), \ldots, S_{1}(e_{4})$ generate $\tau := \cone(e_{1},
  e_{2}, e_{3}) \subset \RR^{3}$. The faces $\cone(e_{1}, e_{3})$ and
  $\cone(e_{2}, e_{3})$ of $\tau$ are not containd in $S_{1}(\vert
  \Delta \vert)$.
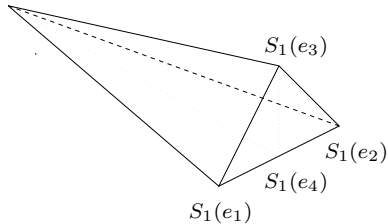
\begin{center}
  \begin{picture}(0,0)%
\includegraphics{nonsurj.pstex}%
\end{picture}%
\setlength{\unitlength}{0.00035000in}%
\begingroup\makeatletter\ifx\SetFigFont\undefined
\def\x#1#2#3#4#5#6#7\relax{\def\x{#1#2#3#4#5#6}}%
\expandafter\x\fmtname xxxxxx\relax \def\y{splain}%
\ifx\x\y   
\gdef\SetFigFont#1#2#3{%
  \ifnum #1<17\tiny\else \ifnum #1<20\small\else
  \ifnum #1<24\normalsize\else \ifnum #1<29\large\else
  \ifnum #1<34\Large\else \ifnum #1<41\LARGE\else
     \huge\fi\fi\fi\fi\fi\fi
  \csname #3\endcsname}%
\else
\gdef\SetFigFont#1#2#3{\begingroup
  \count@#1\relax \ifnum 25<\count@\count@25\fi
  \def\x{\endgroup\@setsize\SetFigFont{#2pt}}%
  \expandafter\x
    \csname \romannumeral\the\count@ pt\expandafter\endcsname
    \csname @\romannumeral\the\count@ pt\endcsname
  \csname #3\endcsname}%
\fi
\fi\endgroup
\begin{picture}(4994,3220)(879,-3259)
\put(3601,-3211){\makebox(0,0)[lb]{\smash{\SetFigFont{8}{9.6}{rm}$S_1(e_1)$}}}
\put(4726,-2761){\makebox(0,0)[lb]{\smash{\SetFigFont{8}{9.6}{rm}$S_1(e_4)$}}}
\put(4726,-736){\makebox(0,0)[lb]{\smash{\SetFigFont{8}{9.6}{rm}$S_1(e_3)$}}}
\put(5626,-2311){\makebox(0,0)[lb]{\smash{\SetFigFont{8}{9.6}{rm}$S_1(e_2)$}}}
\end{picture}

\end{center}
By \cite{acha}, the toric morphism $s_{1} \colon X_{\Delta} \to
X_{\tau}$ defined by $S_{1}$ is the toric quotient for the action of
the subtorus $H \subset T$ corresponding to the sublattice
$\ker(S_{1})$ of $\ZZ^{4}$. In particular, $s_{1}$ is not surjective.
So the action of $H$ on $X_{\Delta}$ has no quotient in the category
of quasi--projective varieties.  \quad $\diamondsuit$
\end{example}

Note that surjectivity of the toric quotient $s_{1} \colon X \to X \tq
H$ does not imply the existence of a quotient in the category of
quasi--projective varieties:

\begin{example}\label{subtil}
  Let $\Delta'$ be the fan in $\RR^{3}$ with the maximal cones
  $$
  \tau_{1} := \cone(e_{1},e_{2}), \quad \tau_{2} :=
  \cone(e_{3},e_{4}), \quad \tau_{3} := \cone(e_{5},e_{6}).$$
  Then the
  associated toric variety $X_{\Delta'}$ is an open toric subvariety
  of $\CC^{6}$. In the notation of Example \ref{nonsurj}, define a
  projection $S_{1} \colon \ZZ^{6} \to \ZZ^{3}$ by
  $$\begin{array}{ll}
    S_{1}(e_{1}) := -e_{1}, & \qquad S_{1}(e_{2}) := v_{1}', \\
    S_{1}(e_{3}) := -e_{2}, & \qquad S_{1}(e_{4}) := v_{2}', \\
    S_{1}(e_{5}) := -e_{3}, & \qquad S_{1}(e_{6}) := v_{3}'. \\
\end{array}$$
Then $S_{1}$ is a map of the fan $\Delta'$ and the fan $\Delta$ of
\ref{nonsurj}, in fact the (surjective) toric morphism $s_{1} \colon
X_{\Delta'} \to X_{\Delta}$ associated to $S_{1}$ is the toric
quotient of the action of the subtorus $H$ of ${\CC^{*}}^{6}$
corresponding to the sublattice $\ker(S_{1}) \subset \ZZ^{6}$. Since
the quasi--projective reduction of $X_{\Delta}$ is not surjective,
there is no quotient in the category of quasi--projective varieties for
the action of $H$ on $X_{\Delta'}$.
\end{example}

\bibliography{}

\bigskip

\bigskip

{\sc
\noindent Fakult\"at f\"ur Mathematik und Informatik \\
\noindent Universit\"at Konstanz \\
\noindent Fach D197 \\
\noindent D-78457 Konstanz \\
\noindent Germany}

\medskip

\noindent{\it E-mail address:}\enspace {\tt
  Annette.ACampo@uni-konstanz.de} \\
\noindent$\hphantom{\hbox{\it E-mail address:}}$\enspace {\tt
  Juergen.Hausen@uni-konstanz.de}

\end{document}